\documentclass[]{article}
\usepackage{vmargin}
\setmarginsrb{3.2cm}{2.5cm}{3.2cm}{2.5cm}{1cm}{1cm}{1cm}{1.6cm}
\usepackage{wasysym}
\usepackage{amsmath}
\usepackage{amssymb}

\newcounter{Prop}
\newcounter{Theo}
\newcounter{Def}

\def\Proposition#1#2{\refstepcounter{Prop}
  \[\left[
  \parbox[l]{15cm}{
    \begin{minipage}[l]{15cm}
      { \sf $\,$ Proposition \theProp:} {\rm #1}
      \begin{center}
        \begin{minipage}[c]{14cm}
          {\em #2}
        \end{minipage}
      \end{center}
    \end{minipage}
    }
\right.\]
}

\def\Proof#1{
  \underline{{\sf Proof:}} 
  {\rm #1}
  \smallskip
  \hfill $\Box$
  }

\def\id{{\rm id}}
\def\End{{\rm End}}
\def\map{{\rm Map}}
\def\deg{{\rm deg}}
\def\cplx{\mathbb{C}}
\def\fA{\mathfrak{A}}  \def\fB{\mathfrak{B}}  
    
\def\fG{\mathfrak{G}}  \def\fH{\mathfrak{H}}  
     
  \def\fN{\mathfrak{N}}  
     
    \def\fU{\mathfrak{U}}

\def\cZ{{\cal{Z}}}

\title{{\bf\sf
    Universal Solutions of Quantum Dynamical\\ Yang--Baxter Equations}}
\author{{\Large
    D. Arnaudon}\thanks{e-mail: daniel.arnaudon@lapp.in2p3.fr}\cr 
  Laboratoire de Physique Th{\'e}orique ENSLAPP\footnote{Laboratoire 
    URA 1436 du CNRS, associated to ENS Lyon and
    Universit{\'e} de Savoie} \cr Chemin de Bellevue,
  BP 110 \cr F-74941 Annecy-le-Vieux Cedex, France \cr \cr
  {\Large E. Buffenoir}\thanks{e-mail: buffenoi@lpm.univ-montp2.fr}\cr 
  Laboratoire de Physique Math{\'e}matique et
    Th{\'e}orique\footnote{Laboratoire UMR 5825 du CNRS} 
  \cr Universit{\'e} Montpellier 2,  Place Eug{\`e}ne Bataillon \cr 
  F-34000 Montpellier France \cr \cr
  {\Large E. Ragoucy}\thanks{e-mail: eric.ragoucy@cern.ch, on leave of
    absence from Laboratoire de Physique Th{\'e}orique ENSLAPP}
  {\ and}
  {\Large Ph. Roche}\thanks{e-mail: philippe.roche@cern.ch, on leave of
    absence from CPT, Ecole Polytechnique F-91128 Palaiseau cedex,
    Laboratoire UPR 14 du CNRS}
  \cr TH Division CERN \cr CH-1211 Geneva 23 Switzerland}
\date{December 12, 1997}

\begin{document}

\maketitle

\begin{abstract}
{\normalsize

\indent

  We construct a universal trigonometric solution of the 
  Gervais--Neveu--Felder equation in the case of finite dimensional 
  simple Lie algebras and finite dimensional contragredient
  simple Lie superalgebras.}
\end{abstract}

\vfil\vfil
\rightline{CERN-TH/97-363}
\rightline{ENSLAPP-A-664/97}
\rightline{q-alg/9712037}

\newpage

\section{Introduction}

Let $\mathfrak{H}$ be a finite dimensional commutative Lie algebra
over 
$\mathbb{C}$, $V$ a semisimple finite
dimensional $\mathfrak{H}-$module. The quantum dynamical Yang--Baxter 
equation also known as 
Gervais--Neveu--Felder
equation (GNF) is:
\begin{equation}
  R_{12}(-\mu+2 h^{(3)})R_{13}(-\mu)R_{23}(-\mu+2h^{(1)})=
  R_{23}(-\mu)R_{13}(-\mu+2h^{(2)})R_{12}(-\mu) \label{GNF}
\end{equation}
where $R:\mathfrak{H}^{*}\rightarrow \End(V\otimes V)$ is a 
function (often chosen to be meromorphic) and where by
definition
\begin{equation}
  R_{12}(-\mu+2h^{(3)})(v_1\otimes v_2\otimes v_3)=
  R_{12}(-\mu+2\eta)(v_1\otimes v_2)\otimes v_3
\end{equation}
if $v_3$ has weight $\eta.$

This equation was first discovered by J.L. Gervais and A. Neveu in
their work on
 the quantisation of Liouville theory
\cite{GN}. It reappeared more recently in 
its modern form in the work of Felder \cite{Fe} in his approach to
 the quantisation of Knizhnik--Zamolodchikov--Bernard equations.
Since then, it has been shown to be one of the basic tools in the
R-matrix formalism of the quantisation of a wide family of models
(Calogero--Moser, Calogero--Sutherland, Ruijsenaars--Schneider) \cite{ABB}.

The aim of our work is to associate to each finite dimensional 
simple Lie algebra (and also to each finite dimensional contragredient
simple Lie superalgebra) a trigonometric universal solution of GNF.

Let $\mathfrak{G}$ be a
simple Lie algebra over $\mathbb{C}$ of rank $r$ and $\mathfrak{H}$ be
a Cartan 
subalgebra of $\mathfrak{G}$. Denote by $(\ell_i)$ an orthonormal
basis of 
$\mathfrak{H}$ with respect to the Killing
form and $\ell^{i}\in \mathfrak{H}^{\star}$ its dual basis. A universal 
solution of GNF is a meromorphic map
$R:\mathbb{C}^r\rightarrow 
\mathfrak{U}_q(\mathfrak{G})^{\otimes 2}$ solving the equation
\begin{equation}
  R_{12}(xq^{2\ell^{(3)}})R_{13}(x)R_{23}(xq^{2\ell^{(1)}})=
  R_{23}(x)R_{13}(xq^{2\ell^{(2)}})R_{12}(x) \label{sYB}
\end{equation}
where $x=(x_1,\cdots,x_r)\in \mathbb{C}^r,$ and
$xq^{\ell}=(x_1q^{\ell_1},\cdots,x_r 
q^{\ell_r})$.
Note that if  $R$ is meromorphic there is no difficulty in defining 
$R_{12}(xq^{2\ell^{(3)}})$.

If $\pi$ is a finite dimensional representation of 
$\mathfrak{U}_q(\mathfrak{G})$ acting on $V$ then $(\pi\otimes
\pi)(R(x))$ is a solution of GNF with the change of variable 
$x_i=q^{-\mu_i}$ where 
$\mu=\sum_{i}\mu_i \ell^i$.

A universal solution of GNF equation can be obtained from a 
solution of the shifted cocycle equation (\ref{eq:s-cocycle}). Let us
assume that 
$F:\mathbb{C}^r\rightarrow \mathfrak{U}_q(\mathfrak{G})^{\otimes 2}$ is an 
invertible element of  weight $0$, i.e. 
$[F_{12},h\otimes 1+1\otimes h]=0,\forall h\in \mathfrak{H},$
satisfying the following
equation:
\begin{equation}
  (\Delta\otimes \id)(F(x))\ F_{12}(xq^{\ell^{(3)}})=
  (\id\otimes \Delta)(F(x))\ F_{23}(x)
  \label{eq:s-cocycle}
\end{equation}
Then $R(x)=F_{21}(x)^{-1}R_{12}F_{12}(x)$ is easily shown to satisfy 
the universal GNF equation where $R$ is a universal
element in $\mathfrak{U}_q(\mathfrak{G})^{\otimes 2}$ satisfying the 
quasi-triangularity axioms.

An explicit formula for $R$  is known for every 
$\mathfrak{U}_q(\mathfrak{G})$ \cite{Dr,Ro,Bu,KR,LS}, see also \cite{KT} 
in the case of quantum superalgebras.
As far as $R(x)$ is concerned, matricial solutions of GNF have been
given in \cite{BG} for $sl(n)$. Elliptic matricial solutions of 
GNF with spectral
parameter have been exhibited in \cite{Fe} and partially 
classified in \cite{EV}.
Explicit solutions to the shifted cocycle are known in the $sl(2)$ 
case \cite{Ba} and also in the $osp(1|2)$ case \cite{BT}.
Unicity of the solution of the shifted cocycle equation 
under some hypothesis and recursion relations have been 
written in \cite{Fr} for general simple Lie algebra.
An  important formula  $F(x)=\prod_{k=0}^{+\infty}F_k$ has been obtained,
and $F_k$ has been exactly computed in the case of $sl(2)$. 
In the rest of this work we will give exact formulas for $F_k$ and provide
an alternative approach to the computation of $F(x)$ (using a linear
equation) and to the proof of the
cocycle identity (performing an algebraic as well as an analytical study).

\section{Notations}

In the following we will always assume that $q\in \mathbb{R}$, with $0<q<1.$
We will denote by $\kappa$ the restriction of the Killing 
form to $\mathfrak{H}$. 
If $\alpha\in\mathfrak{H}^{\star}$
we will denote by $t_{\alpha}\in \mathfrak{H}$ the element defined by 
$\kappa(t_{\alpha},h)=\alpha(h),\forall h\in \mathfrak{H},$ and by 
$(.\vert.)$ the scalar product on $\mathfrak{H}^{\star}$
defined from $\kappa$ by duality.
Let $(\alpha_{i},i=1,\dots,r)$ be a choice of simple roots,
$\Phi$ the set of roots and $\Phi^+$ the
corresponding set of positive roots. 
To each root $\alpha$ we will associate the element 
$h_{\alpha}=\frac{2}{(\alpha\vert \alpha)}t_{\alpha}.$ 

A presentation of $\mathfrak{U}_q(\mathfrak{G})$ by generators and
relations is given by:
\begin{eqnarray}
  {[t_{\alpha_i},t_{\alpha_j}]}\  =\  0 &&
  {[e_{\alpha_i},f_{\alpha_j}]}\  =\  \delta_{ij}\  
  \frac{q^{t_{\alpha_i}}-q^{-t_{\alpha_i}}}{q-q^{-1}}
\label{defUq-1}\\
  \begin{array}{l}
  {[t_{\alpha_i},e_{\alpha_j}]} \ =\  a^{\mbox{\tiny sym}}_{ij}
  e_{\alpha_j}\\
  {[t_{\alpha_i},f_{\alpha_j}]}\ =\ -a^{\mbox{\tiny sym}}_{ij} f_{\alpha_j} \\
  \end{array}
  &\mbox{ with }& a^{\mbox{\tiny sym}}_{ij}=(\alpha_i|\alpha_j)\\
  \begin{array}{l}
    ({\rm ad}_{q'} e_{\alpha_i})^{n_{ij}}(e_{\alpha_j})\ =\ 0\\[.21cm]
    ({\rm ad}_{q'} f_{\alpha_i})^{n_{ij}}(f_{\alpha_j})\ =\ 0
  \end{array}
  & \mbox{ if }& i\neq j \mbox{ and } 
  n_{ij}=1-2\frac{a^{\mbox{\tiny sym}}_{ij}}{a^{\mbox{\tiny sym}}_{ii}},\ 
  q'=q \mbox{ or } q^{-1}\label{defUq-2}
\end{eqnarray}
where we have introduced:
\begin{equation}
({\rm ad}_{q^{\pm1}} x) (y)= \sum_{(x)} x_{(1)}\, y\, S^{\pm1}(x_{(2)}).
\end{equation}
We will call $\fU_q(\fN_+)$ (resp. $\fU_q(\fN_-)$) the subalgebra of 
$\mathfrak{U}_q(\mathfrak{G})$ generated
by $e_{\alpha_i},\ i=1,\dots,r$ (resp. $f_{\alpha_i},\ i=1,\dots,r$).
As usual, we denote by $\mathfrak{U}_q(\fB_+)$ 
(resp. $\mathfrak{U}_q(\fB_-)$) the
algebra generated by $h_{\alpha_i}$, $e_{\alpha_i}$ $i=1,\dots,r$
(resp. $h_{\alpha_i}$, $f_{\alpha_i}$ $i=1,\dots,r$).

The Hopf algebra structure is defined by:
\begin{eqnarray}
  \Delta(t_{\alpha_i}) &=& t_{\alpha_i}\otimes 1+1\otimes t_{\alpha_i}\\
  \Delta(e_{\alpha_i})&=& e_{\alpha_i}\otimes q^{t_{\alpha_i}}
  +1\otimes e_{\alpha_i}\\
  \Delta(f_{\alpha_i}) &=& f_{\alpha_i}\otimes 1
  +q^{-t_{\alpha_i}}\otimes f_{\alpha_i}
\end{eqnarray}

Let us fix a
normal ordering $<$ on the set of positive roots (see \cite{KT} 
and references therein). Such a normal ordering is defined by
requiring that for $\alpha,\ \beta,\ \alpha+\beta$ positive roots
such that $\alpha<\beta$, we have $\alpha<\alpha+\beta<\beta$.
For a given set of positive roots, there exists several normal
orderings. It can be shown that they are related by elementary moves
recalled in \cite{KT}.
Then, to each non-simple positive root $\alpha$ we can
associate elements (which depend on the choice of the normal
ordering) $e_\alpha\in\mathfrak{U}_q(\mathfrak{N}_+)$ and
$f_\alpha\in\mathfrak{U}_q(\mathfrak{N}_-)$. They are uniquely defined
by:
\begin{equation}
  e_{\alpha+\beta} \ =\  e_\alpha e_\beta -q^{-(\alpha|\beta)} e_\beta
  e_\alpha \mbox{ and }
  f_{\alpha+\beta} \ =\ f_\beta f_\alpha -q^{(\alpha|\beta)} 
  f_\alpha f_\beta
\end{equation}
for $\alpha<\beta$ such that there is no pair
$\{\alpha',\beta'\}$ in $]\alpha,\beta[$ with
$\alpha+\beta=\alpha'+\beta'$.
It can be shown that 
\begin{equation}
  [e_\alpha,f_\alpha]=a_\alpha\ \frac{q^{t_\alpha}-q^{-t_\alpha}}{q-q^{-1}}
\end{equation}
where $a_\alpha$ are functions of $q$.

A Poincar{\'e}--Birkhoff--Witt basis of $\mathfrak{U}_q(\mathfrak{N}_+)$
(resp. $\mathfrak{U}_q(\mathfrak{N}_-)$) is given by:
$(e^p=\prod_{\alpha\in\Phi^+}^{>} (e_\alpha)^{p_\alpha})_{p\in\cZ}$
(resp. $(f^{r}=\prod_{\alpha\in\Phi^+}^{>}
(f_\alpha)^{r_\alpha})_{r\in\cZ}$), 
where $\cZ=\map(\Phi^+,\mathbb{Z}^+)$
and the products are ordered according to $>$, the reversed normal order
on $\Phi_+$.

As a notation, for $p\in\cZ$, we introduce 
$\rho_p=\frac{1}{2}\sum_{\alpha\in\Phi^+}\ p_\alpha \alpha$ which
satisfies
\begin{equation}
  {}[h,e^p]=\rho_p(h)\ e^p \mbox{ and } [h,f^p]=-\rho_p(h)\ f^p,\ \
  \forall h\in\mathfrak{H}.
\end{equation}

As usual, we denote by
$\rho=\frac{1}{2}\sum_{\alpha\in\Phi^+}\alpha$ the Weyl vector.
$\mathfrak{U}_q(\mathfrak{G})$ is $\mathbb{Z}$-graded with respect 
to the adjoint action of $t_\rho$, 
which is called the principal gradation of
$\mathfrak{U}_q(\mathfrak{G})$. 
The principal degree of the generators $e_{\alpha_i}$ is 1.

$\mathfrak{U}_q(\mathfrak{G})$ is quasi-triangular and let 
$R\in\left(\fU_q(\fB_+)\otimes\fU_q(\fB_-)\right)^c$ be a 
solution of the quasi-triangularity
conditions.
We have denoted with a supercript $c$ a completion of 
$\fU_q(\fB_+)\otimes\fU_q(\fB_-)$ in the
usual sense: elements of $(\fU_q(\fB_+)\otimes\fU_q(\fB_-))^c$ 
are of the form $g=\sum_{p,r\in\cZ} (e_p\otimes f_r) \psi_{p,r}(x)$
where $\psi_{p,r}(x)$ are arbitrary functions of $h_{(1)}$
and $h_{(2)}$, i.e. are elements
of \map$(\fH^*\oplus\fH^*,\cplx)$.
It can be shown that we have the explicit multiplicative formula 
\begin{equation}
R=K{\hat R} \ \mbox{ where }\ 
K=\prod_{j=1}^r q^{\ell_j\otimes\ell_j}\ \mbox{ and }\
{\hat R}=\prod_{\alpha\in \Phi^+}^{>}{\hat R}_{\alpha}
\label{R=KR}
\end{equation}
and ${\hat R}_{\alpha}=\exp_{{\bar q}_{\alpha}}((q-q^{-1})
a_{\alpha}^{-1}e_{\alpha}\otimes f_{\alpha})$
with ${\bar q}_{\alpha}=q_\alpha^{-(\alpha|\alpha)}$.
The $q$-exponential is defined by:
\begin{equation}
\exp_q(z)\ =\ \sum_{n=0}^{+\infty} \frac{1}{[n]_q!}\ z^n
\mbox{ with } [n]_q!=[n]_q\, [n-1]_q \cdots [2]_q [1]_q
\mbox{ and }[n]_q=\frac{1-q^{n}}{1-q}.
\end{equation}
An important result states that the value of $R$ is independent of
the choice of normal ordering.

{\it Remark:} We have essentially used the conventions of \cite{KT},
except that we have applied to all the expressions therein the
antimorphism and comorphism $\dagger$ defined by
\begin{equation}
  \underline{e}_\alpha^\dagger = f_\alpha \mbox{ , }  
  \underline{f}_\alpha^\dagger = e_\alpha 
  \mbox{ and } \underline{h}_\alpha^\dagger = t_\alpha \ 
  \forall \alpha\in\Phi_+
\end{equation}
where the underlined expressions denote the bases used in \cite{KT}.
The $R$ matrix we use is 
$R_{12}=({\dagger\otimes\dagger})(\underline{R}_{21})$.

\section{Universal solution in the case of Lie algebras}

The aim of our work is to show that, hidden in the resolution of the
shifted cocycle equation, 
there is a linear equation satisfied by $F$, 
the solution of which is unique
under some hypothesis. This linear equation was shown to be satisfied by
$F(x)$ in the $sl(2)$ case in \cite{BR}. The necessity of this
equation 
came from a complete different point of view: in 
\cite{BR} it
has been shown that the characters of the irreducible unitary 
representations of the quantum Lorentz
group are constructed from complex continuation of $6j$ symbols of 
$\mathfrak{U}_q(su(2)).$ It has been shown in 
\cite{BBB} that there is a complete dictionary between matrix
elements 
of $F(x)$ and complex continuation of $6j$ symbols.
{}From the regularity conditions of the characters of the quantum
Lorentz 
group and the explicit expression 
of $F$,
a linear equation satisfied by $F$ was derived.

Let $B(x)\in \fU_q(\fH)^c={\rm \map}(\mathfrak{H}^*,\cplx)$ be the element 
$B(x)=\prod_{j=1}^r (x_{j}^{\ell_j}q^{\ell_j^2})=
q^{\sum_{j=1}^r(\ell_j\ell_j-\mu_j\ell_j)},$  
the generalisation of the linear equation of  \cite{BR} can be written:
\begin{equation}
  F_{12}(x)B_{2}(x)={\hat R}_{12}^{-1}B_{2}(x)F_{12}(x)\label{lineareq}
\end{equation}

We first have the following proposition on the existence and unicity
of 
the solution of this linear equation.

\Proposition{}{\label{Ftriang}
  Under the following hypothesis:
  \begin{enumerate}
  \item $F(x)\in \left(\mathfrak{U}_q(\mathfrak{B}^+)\otimes 
    \mathfrak{U}_q(\mathfrak{B}^-)\right)^c$
  \item the projection of $F(x)$ on
    $\left(\mathfrak{U}_q(\mathfrak{H})^{\otimes 2}\right)^c$ 
    is equal to $1\otimes 1$
  \end{enumerate}
  there exists a unique solution of equation (\ref{lineareq}).\\
  Moreover, under these assumptions, we have 
  \begin{equation}
    F(x)=\sum_{p,r\in\cZ } e^{p}\otimes \left(f^{r}\ \phi_{p,r}(x)\right)
    \label{serie_som}
  \end{equation} 
  where  $\phi_{p,r}(x)$ are
  rational functions on $\mathfrak{H}$.\\
  Under these conditions $F$ is of weight zero i.e.
  $[F(x), h\otimes 1+1\otimes h]=0,\ 
  \forall h\in \mathfrak{H}.$
  }
\Proof{
  Starting with ${\hat R}_{12}^{-1}=\sum_{s\in\cZ}\sigma_s\ e^s\otimes
  f^s$ and 
$F(x)=\sum_{p,r\in\cZ }\ (e^{p}\otimes f^{r})\ \psi_{p,r}(x)$, where  
  $\psi_{p,r}(x)$ 
  a priori belongs to  $\map(\mathfrak{H}^*\oplus
  \mathfrak{H}^*,\cplx)$, 
  one computes that the equation (\ref{lineareq}) can be rewritten as
  \begin{equation}
    \left(1-q^{(\mu+2\rho_{r}|2\rho_{r})}q^{-4t^{(2)}_{\rho_{r}}}
    \right)
    \psi_{p,r}(x)=
    \sum_{\begin{subarray}{c}k+s=p\\ l+s=r\\ s\neq0\end{subarray}} 
    a_{p}^{ks} b_{r}^{ls} \sigma_s 
    q^{(\mu+2\rho_{l}|2\rho_{l})}q^{-4t^{(2)}_{\rho_{l}}}
    \psi_{k,l}(x)
    \label{recurence}
  \end{equation}
  where we have introduced the numbers $a_{p}^{ks}$ and $b_{r}^{ls}$
  defined by  
  \begin{equation*}
    e^k e^s \ =\ \sum_{p\in\cZ} a_{p}^{ks} e^p \ \mbox{ and }\ 
    f^k f^s \ =\ \sum_{p\in\cZ} b_{p}^{ks} f^p
  \end{equation*}
  The equation (\ref{recurence}) clearly shows that we have a linear 
  system which is strictly triangular 
  with respect to the principal gradation. Moreover, once we have 
  normalised $\psi_{0,0}(x)=1$, 
  the equation proves that $\psi_{p,r}$ are rational functions on 
  $1\otimes\mathfrak{H}$ 
  which are uniquely defined.\\
  Now, if $F_{12}(x)$ is a solution of (\ref{lineareq}) satisfying 
  the two hypothesis of proposition 
  \ref{Ftriang}, 
  for any $\lambda\in\cplx$ and any element 
  $h\in\mathfrak{H}$, 
  $\tilde{F}_{12}(x)=
  \lambda^{h_{(1)}+h_{(2)}}F_{12}(x)\lambda^{-(h_{(1)}+h_{(2)})}$ 
  is also a solution of 
  (\ref{lineareq}) satisfying the two hypothesis of proposition
  \ref{Ftriang}.  
  Thus, as the solution is unique, we must have
  $\tilde{F}_{12}(x)=F_{12}(x)$, 
  so that $F$ is of zero weight.
  }  

Although we have not mentioned up to now the problems related to the 
convergence of the series 
(\ref{serie_som}), it suffices to say that in
each finite dimensional representation this series has a finite number
of terms. 
In the sequel when we will write a series or an infinite product we
will always assume that they are evaluated in finite
dimensional representation, and we will therefore only have to
consider the convergence of a series or infinite product 
of finite dimensional matrices.

We now show that the unique solution of (\ref{lineareq})
satisfying the relations (\ref{eq:s-cocycle}) 
can be written as an infinite and very explicit product.
Let us define the formal product 
\begin{equation}
  F=\prod_{k=0}^{+\infty}B_{2}^k {\hat R}_{12}^{-1}B_{2}^{-k}
  \label{Fprod}
\end{equation} 
the convergence of which is analysed
in the  proposition \ref{Fconv}. We will show that $F(x)$ satisfies the
linear equation (\ref{lineareq}) and moreover satisfies 
the cocycle identity (\ref{eq:s-cocycle}).
This proposition is divided in two parts: an algebraic study and an
analytic 
one where the limit $N\rightarrow +\infty$ is analysed.

Let us define, for each $N\in \mathbb{N},\,\,{\buildrel N\over F}(x)=
\prod_{k=0}^{N}B_{2}^k(x) {\hat R}_{12}^{-1}B_{2}^{-k}(x)$.
\Proposition{}{
  ${\buildrel N\over F}(x)$ is the unique sequence of elements of 
  $\mathfrak{U}_q(\mathfrak{G})^{\otimes 2}$ satisfying the recursion
  equation 
  \begin{equation}
    {\buildrel N\over F}(x)B_2(x)={\hat R}_{12}^{-1}B_2(x){\buildrel
    N-1\over F}(x) \label{Frecursion}
  \end{equation}
  with ${\buildrel 0\over F}(x)={\hat R}_{12}^{-1}.$
  }

\Proof{Trivial computation.}

Remark: we will show in the proposition \ref{Fconv} that, evaluated in the
tensor product of any finite dimensional representation, the
product (\ref{Fprod}) is convergent for $\mu$ dominant weight 
with all the scalar 
products $(\mu\vert \alpha_i)$ 
sufficiently large. As a result, in this domain, one can show from
(\ref{Frecursion}) that 
$\lim_{N\rightarrow+\infty} {\buildrel N\over  F}(x)$ satisfies 
the linear equation, and therefore is
equal to the expression (\ref{serie_som}).

\Proposition{}{
  Let us define $U=B_1B_2K_{12}^2$ and $X=B_2K_{12}^2$.
  We have the following two identities:
  \begin{eqnarray} 
    &&{\buildrel N\over F}_{23}(x){\buildrel N\over
      F}_{12}{}^{-1}(xq^{\ell^{(3)}})
    \ =\ 
    \left(\prod_{k=0}^{N} B_{3}^k {R}_{23}^{-1} {K}_{23} B_{3}^{-k}
    \right)
    \left(\prod_{k=N}^0 U_{23}^{k} K_{12}^{-1} {R}_{12} U_{23}^{-k}
    \right)\\
    &&(\id\otimes\Delta)({\buildrel N\over F}{}^{-1}(x))\ 
    (\Delta\otimes\id)({\buildrel N\over F}(x))\ =\ 
\label{eq24}\\
    &&\ \left(\prod_{k=0}^{N}B_{3}^k R_{23}^{-1}(X_{23}^{N+1} 
    {\hat R}_{13}^{-1}X_{23}^{-N-1})
    K_{23}B_{3}^{-k}\right)
    \left(\prod_{k=N}^0 U_{23}^{k}K_{12}^{-1}(B_{3}^{N+1}
      {\hat R}_{13}B_{3}^{-N-1})R_{12}U_{23}^{-k}
    \right)
    \nonumber
  \end{eqnarray}
  }

\Proof{
  The first relation is an immediate application of the relation
  $B_{2}(xq^{\ell^{(3)}})=B_{2}(x)K_{23}^2.$ 
  The second relation is more tricky.
  Let us define $V=K_{12}^{-1}K_{13}^{-1}R_{13}R_{12}$ and 
  $W=R_{23}^{-1}R_{13}^{-1}K_{13}K_{23}B_3,$
  we have the following identity:
  \begin{equation}
    U_{23}^{-1}VW=WU_{23}^{-1}V.\label{UVW}
  \end{equation}
  This identity is shown to be equivalent to Yang--Baxter equation
  with  the use of the two relations: 
  \begin{equation}
    R_{12}K_{13}K_{23}=K_{13}K_{23}R_{12} 
    \mbox{ and } R_{12}U_{12}=U_{12}R_{12}
  \end{equation}
  From the definition of ${\buildrel N\over F}$ we have 
  \begin{eqnarray}
    (\id\otimes \Delta)({\buildrel N\over F}{}^{-1}) &=&
    \prod_{k=N}^0 U_{23}^k V U_{23}^{-k}=U_{23}^{N+1}(U_{23}^{-1}V)^{N+1}\\
    (\Delta\otimes\id)({\buildrel N\over F}) &=&
    \prod_{k=0}^{N}B_{3}^k WB_{3}^{-k-1}=W^{N+1}B_{3}^{-N-1}.
  \end{eqnarray}
  We can therefore write, using the identity (\ref{UVW}):
  \begin{eqnarray}
    &&(\id\otimes \Delta)({\buildrel N\over F}{}^{-1})\ 
    (\Delta\otimes\id)({\buildrel N\over F})
    \ =\ U_{23}^{N+1}W^{N+1}(U_{23}^{-1}V)^{N+1}B_{3}^{-N-1}\\
    &&\ \ =(U_{23}^{N+1}WU_{23}^{-N-1})^{N+1}\, U_{23}^{N+1}\, 
    (U_{23}^{-1}V)^{N+1}\, B_{3}^{-N-1}\\
    &&\ \ =\left(R_{23}^{-1}X_{23}^{N+1}{\hat R}_{13}{}^{-1}
    X_{23}^{-N-1}K_{23}B_{3}\right)^{N+1}\, 
    U_{23}^{N+1}\, (U_{23}^{-1}V)^{N+1}\, B_{3}^{-N-1}\\
    &&\ \ =\left(\prod_{k=0}^{N}B_{3}^k R_{23}^{-1}(X_{23}^{N+1} 
    {\hat R}_{13}^{-1}X_{23}^{-N-1}
    )K_{23}B_{3}^{-k}\right)B_{3}^{N+1}U_{23}^{N+1}(U_{23}^{-1}V)^{N+1}
    B_{3}^{-N-1}\\
    &&\ \ =\left(\prod_{k=0}^{N}B_{3}^k R_{23}^{-1}(X_{23}^{N+1} 
    {\hat R}_{13}^{-1}X_{23}^{-N-1}
    )K_{23}B_{3}^{-k}\right)\left(\prod_{k=N}^{0}U_{23}^{k}B_{3}^{N+1}
    VB_{3}^{-N-1}U_{23}^{-k}\right)
  \end{eqnarray}
  which ends the proof of this proposition.
  }

{}From this result, in order to prove that $\lim_{N\rightarrow+\infty}
{\buildrel N\over  F}(x)$  
satisfies the cocycle
identity, it is sufficient to show that both 
$B_{3}^{N+1}{\hat R}_{13}B_{3}^{-N-1}$  and  
$X_{23}^{N+1}{\hat R}_{13}{}^{-1}X_{23}^{-N-1}$ tend to $1$
sufficiently fast. This is indeed the case for the solution (\ref{Fprod}) 
in each finite dimensional representation of
$\mathfrak{U}_q(\mathfrak{G}).$ 

We will need the following proposition:

\Proposition{}{\label{banach}
  \begin{enumerate}
  \item 
    Let $\fA$ be a Banach algebra with norm $\|.\|$, and let
    $(u_{k})_{k\in\mathbb{N}}$ be a sequence of elements 
    of $\fA.$ A sufficient condition for the product
    $\prod_{k=0}^{+\infty}u_k$ to be convergent is: 
    \begin{itemize}
    \item $\exists C>0, \forall n>0, \forall m\leq n,\ 
     \sum_{k=m}^{n}\log(\|u_k\|)\leq C$.
    \item $\sum_{k=0}^{+\infty}\|u_k-1\|$ is convergent. 
    \end{itemize}
  \item  
    Let
    $(v_{k}^{(n)})_{(k,n)\in\mathbb{N}^2}$ be a  
    sequence of elements of $\fA$. The product
    $\prod_{k=0}^{n}u_k v_{k}^{(n)}$ converges to
    $\prod_{k=0}^{+\infty}u_k$ if the previous assumptions are  
    satisfied, together with 
    \begin{itemize}
    \item $\exists C>0, \forall n>0, \forall m\leq n,\ \prod_{k=m}^{n}
      \|v_{k}^{(n)}\|\leq C$\\ 
    \item $\exists C'>0, \forall n>0, \forall k\leq n,\
      \|v_{k}^{(n)}-1\|\leq \frac{C'}{n^2}$. 
    \end{itemize}
  \end{enumerate}
 These sufficient assumptions are of course not at all minimal but we
 will only need these crude hypothesis. 
  }

\Proof{
  The first part of this 
  proposition is a direct application of Cauchy criterion
  to the partial  
  product $\prod_{k=0}^{n}u_k.$ 
  The second one is  a direct application of the inequality:
  \begin{equation}
    \|\prod_{k=0}^{n}u_kv_{k}^{(n)}-\prod_{k=0}^{n}u_k\|\leq 
    \prod_{l=0}^n \|u_l\| \sum_{k=0}^{n}\|v_{k}^{(n)}-1\|\ 
    \|v_{k+1}^{(n)}\|\cdots \|v_{n}^{(n)}\|
  \end{equation}}

We will now apply this proposition to the case of the sequence
$u_{k}=B_{2}^k(x) {\hat R}_{12}^{-1}B_{2}^{-k}(x)$

\Proposition{}{\label{Fconv}
  The product $\prod_{k=0}^{+\infty}B_{2}^k(x) 
  {\hat R}_{12}^{-1}B_{2}^{-k}(x)$ in each finite dimensional
  representation of  
  $\mathfrak{U}_q({\mathfrak{G}})$ is
  convergent for $\mu$ such that all the scalar products
  $(\mu|\alpha_i)$ are sufficiently large.} 

\Proof{
Let $\pi_1$ and $\pi_2$ be finite dimensional representations of
$\fU_q(\fG)$. 
  By iterative application of the proposition \ref{banach} to 
  $u_k=\prod_{\alpha\in\Phi^+}^> B_2^k \hat{R}^{-1}_\alpha B_2^{-k}$,
  it is sufficient to apply the  
  first criterion of proposition \ref{banach} to the sequences 
  $u_k^{\alpha}=(\pi_1\otimes\pi_2)
\left(B_2^k \hat{R}^{-1}_\alpha B_2^{-k}\right)$.\\
  Then, it is easy to compute that 
  \[
  u_k^{\alpha}=\exp_{{\bar q}_\alpha} \left(-a^{-1}_\alpha (q-q^{-1})
  q^{-k(\alpha|\alpha-\mu)} 
  (\pi_1\otimes\pi_2)
\left(e_\alpha\otimes q^{-2kt_\alpha}f_\alpha\right)\right)
  \]
  For $q\in{]0,1[}$, 
  $\|u_k^{\alpha}\|\leq \exp_{{\bar q}_\alpha} 
  \left(\epsilon_\alpha
    q^{k\{(\mu-\alpha|\alpha)-2\|\pi_2(t_\alpha)\|\}} 
  \right)$, with 
  $\epsilon_\alpha=\left| a^{-1}_\alpha (q-q^{-1})\right| 
  \|(\pi_1\otimes\pi_2)(e_\alpha \otimes f_\alpha) \| 
  $. 
  Therefore,
  \begin{equation}
    \log\left(\prod_{k=m}^{k=n}\|u_k^{\alpha}\|\right) \leq
    \sum_{k=m}^{k=n} w_k^\alpha\mbox{ with }w_k^\alpha=
    \log\left(\exp_{{\bar q}_\alpha} 
    \left(\epsilon_\alpha
      q^{k\{(\mu-\alpha|\alpha)-2\|\pi_2(t_\alpha)\|\}}
    \right)\right)\label{prtsum}
  \end{equation}
  Now, if $(\mu|\alpha)> (\alpha|\alpha)+2\|\pi_2(t_\alpha)\|$, we have
  $w_k^\alpha\, 
  \raisebox{-1ex}{${\buildrel \sim \over
      {\scriptstyle{k\rightarrow+\infty}} }$}\, 
  \epsilon_\alpha q^{k\{(\mu-\alpha|\alpha)-2\|\pi_2(t_\alpha)\|\}}$. 
Thus, $\sum_k
  w_k^\alpha$ is convergent,  
  which proves 
  from the inequality (\ref{prtsum}), that the partial sums
  $\sum_{k=m}^n \log(\|u_k^\alpha\|)$ are bounded. 
  \\
  Moreover, 
  $\|u_k^{\alpha}-1 \| \leq 
  \left(\exp_{{\bar q}_\alpha} 
    \left(\epsilon_\alpha
      q^{k\{(\mu-\alpha|\alpha)-2\|\pi_2(t_\alpha)\|\}} 
    \right) -1
  \right)
  \raisebox{-1ex}{${\buildrel \sim \over
      {\scriptstyle{k\rightarrow+\infty}} }$}\, 
  \epsilon_\alpha q^{k\{(\mu-\alpha|\alpha)-2\|\pi_2(t_\alpha)\|\}}$. 
As a result, 
  $\sum_{k=0}^{+\infty} \|u_k^\alpha -1 \|$ is convergent. 
  }

\Proposition{}{
$F$ as defined by (\ref{Fprod})
 satisfies the shifted cocycle equation (\ref{eq:s-cocycle})  
in each finite dimensional representation of  
  $\mathfrak{U}_q({\mathfrak{G}})$.
} 

\Proof{
The two products appearing in the right-hand-side of equation (\ref{eq24})
can be rewritten as 
\begin{eqnarray}
\prod_{k=0}^{N}B_{3}^k{\hat R}_{23}^{-1}X_{23}^{N+1} {\hat R}_{13}^{-1}
    X_{23}^{-N-1}K_{23}B_{3}^{-k} &=&
\prod_{k=0}^{N}B_{3}^k R_{23}^{-1}K_{23}B_{3}^{-k}\ v_{k,N}\\
\prod_{k=N}^{0}U_{23}^{k}K_{12}^{-1}B_{3}^{N+1}
    {\hat R}_{13}B_{3}^{-N-1}R_{12}U_{23}^{-k} &=&
\prod_{k=N}^{0}U_{23}^{k}K_{12}^{-1}R_{12}U_{23}^{-k}\ w_{k,N}
\end{eqnarray}
where $v_{k,N}=\prod^>_{\alpha\in\Phi_+}v^\alpha_{k,N}$ and 
$w_{k,N}=\prod^>_{\alpha\in\Phi_+}w^\alpha_{k,N}$ with 
\begin{eqnarray*}
v_{k,N}^\alpha &=& K^{-1}_{23}B_{3}^k(X_{23}^{N+1} 
    ({\hat R}_{\alpha})_{13}^{-1}X_{23}^{-N-1})K_{23}B_{3}^{-k}\\
w_{k,N}^\alpha &=& U_{23}^{k}K_{12}^{-1}B_{3}^{N+1}K_{13}^{-1}
    ({\hat R}_\alpha)_{13}B_{3}^{-N-1}K_{12}U_{23}^{-k}.
\end{eqnarray*}
In order to prove the above proposition, it is sufficient to show that
the two sequences $v_{k,N}$ and $w_{k,N}$ obey to the criterions
of proposition \ref{banach}.2 for $\mu$ such that all the scalar products
  $(\mu|\alpha_i)$ are sufficiently large.\\
As in proposition \ref{Fconv}, it is sufficient to show that  the sequences 
$v^\alpha_{k,N}$ and $w^\alpha_{k,N}$ satisfy this criterion. 
An easy computation leads to
\begin{eqnarray*}
v_{k,N}^\alpha &=& \exp_{{\bar q}_\alpha} \left(-a^{-1}_\alpha (q-q^{-1})
  q^{-(N+1+k)(\alpha|\alpha-\mu)} 
  q^{-2(N+1+k)t^{(3)}_\alpha-(2N+1)t^{(2)}_\alpha}
  e_\alpha\otimes1\otimes f_\alpha\right)\\
w_{k,N}^\alpha &=& \exp_{{\bar q}_\alpha} \left(a^{-1}_\alpha (q-q^{-1})
  q^{-(N+1+k)(\alpha|\alpha-\mu)} 
  q^{-2(N+1+k)t^{(3)}_\alpha-(2k+1)t^{(2)}_\alpha}
  e_\alpha\otimes 1\otimes f_\alpha\right)
\end{eqnarray*}
For each finite dimensional representation $\pi$, let us denote
$K_\pi=\|q^{-\pi(t_\alpha)}\|$. If $\pi_1$, $\pi_2$ and
$\pi_3$ are finite dimensional representations of $\fU_q(\fG)$
\[
\|\pi_{123}(v_{k,N}^\alpha)-1\|\leq C_{123}\
q^{-(N+k+1)(\alpha|\alpha-\mu)}\ 
(K_{\pi_1})^{2(N+k+1)}\ (K_{\pi_3})^{2k+1} 
\mbox{ with } \pi_{123}=\pi_1\otimes\pi_2\otimes\pi_3
\]
Thus, for $\mu$ such that all the scalar products
  $(\mu|\alpha_i)$ are sufficiently large, there exists 
 $r\in{}]0,1[$ and $K_0>0$ such that 
$\|\pi_{123}(v_{k,N}^\alpha)-1\|\leq K_0\ r^N$, 
which proves that $v_{k,N}^\alpha$
obeys to the second hypothesis of proposition \ref{banach}.2. Now, from 
$\|\pi_{123}(v_{k,N}^\alpha)\|\leq
1+\|\pi_{123}(v_{k,N}^\alpha)-1\|\leq 1+K_0\ r^N$, one
  easily shows that $v_{k,N}^\alpha$ also satisfies to the first
  condition. This ends the demonstration for the sequence $v_{k,N}^\alpha$.
For $w_{k,N}^\alpha$, the proof follows the same steps.
}

\Proposition{}{
A universal solution of GNF is given by
\[
R(x)\ = \ \left(\prod_{k=+\infty}^{0} B_{1}^k {\hat
    R}_{21}B_{1}^{-k}\right)
 R_{12} \left(\prod_{k=0}^{+\infty} B_{2}^k {\hat
    R}_{12}^{-1}B_{2}^{-k}\right).
\]
This solution satisfies 
\[
R_{12}(x)B_2(x)R_{21}(x)=B_2(x)K_{12}^2.
\]}

\Proof{
Trivial computation using the linear equation satisfied by $F$.\\
}

\indent

{\it Example: case of $sl(2)$}

\indent

In \cite{Ba}, a solution to the shifted cocycle equation was
constructed. It reads
\begin{equation}
F_{12}(x)\ =\ \sum_{n=0}^{+\infty}
\frac{(q-q^{-1})^n}{[n]_{q^{2}}!}\, e^n\otimes f^n
\frac{(-1)^n}{\prod_{\nu=1}^{n}\left(1-x^{-2}q^{2\nu}q^{-2h_{(2)}}\right)} 
.
\end{equation}
This solution satisfies the linear equation \cite{BR} and therefore
can also be written as, for $x$ sufficiently large, 
\begin{equation}
F(x) = \prod_{k=0}^{+\infty} B_2(x)^k {\hat R}_{12}^{-1} B_2(x)^{-k}
 = \prod_{k=0}^{+\infty}\ \exp_{q^{2}}\left( -(q-q^{-1})\ x^{-2k}
  q^{-2k} q^{-2kh_{(2)}}\ e\otimes f\right).
\end{equation}
This last formula appeared in 
\cite{Fr}\footnote{However, contrary
  to what has been stated in this work, it cannot be
written as $\exp_{q^{-2}}\left( \phi(x,h_{(2)})\ e\otimes f\right)$.}.

\section{Universal solution for contragredient Lie superalgebras}

The above construction can be applied to the case of 
finite dimensional Lie superalgebras
with the following modifications. We will only consider 
contragredient Lie superalgebras, i.e. 
finite dimensional classical (simple) Lie superalgebras which 
admits a unique non
degenerate invariant supersymmetric bilinear form $(.|.)$ \cite{Ka}. These are
formed by the $A(m,n)$, $B(m,n)$, $C(n+1)$ and $D(m,n)$ infinite series and the
exceptional  $D(2,1;\alpha)$, $G(3)$ and $F(4)$ superalgebras.

$\fU_q(\fG)$ is a
$\mathbb{Z}_2$-graded algebra which admits a presentation by
generators and relations, as in (\ref{defUq-1}-\ref{defUq-2}) 
with the following alterations.
Commutators become graded commutators. The $q$-adjunction is 
also graded. There are supplementary
relations for some types of Dynkin diagrams\cite{SchFLV}. The Hopf
superalgebra structure differs from the usual Hopf algebra structure
by the supercommutativity of the tensor product:
\begin{equation}
(a\otimes b)= (a\otimes 1)(1\otimes b)
=(-1)^{{\rm deg}a\,{\rm deg}b}(1\otimes b)(a\otimes 1)
\end{equation}
where $\deg x$ is the $\mathbb{Z}_2$ degree of the homogeneous element
$x$.
Let $\ell_i$ be an orthogonal basis of $\fH$, the Cartan subalgebra of
$\fU_q(\fG)$. We denote by $\eta$ the restriction of $(.|.)$ to $\fH$:
it is non degenerate and of signature $(p,p')$.
The  case $A(n,n)$ deserves a special treatment, related to the center
of $sl(n|n)$ \cite{KT}: for convenience, we will exclude this case in
the sequel.

With these modifications, it is easy to check from \cite{KT} that we
still have the multiplicative formula (\ref{R=KR}), with now
\begin{eqnarray}
K &=& q^{\sum_{i,j}\eta^{ij}\ell_i\otimes\ell_j}
\ \mbox{ where }\ \eta^{ij}=(\eta^{-1})_{ij} \\
{\hat R}_\alpha &=& \exp_{{\bar q}_\alpha}\left( 
(-1)^{\deg\alpha} (q-q^{-1}) e_\alpha\otimes f_\alpha \right)
\ \mbox{ with }\  {\bar q}_\alpha= (-1)^{\deg\alpha}q^{-(\alpha|\alpha)}.
\end{eqnarray}
The linear equation we consider is identical to (\ref{lineareq}) provided one
uses $B(x)=q^{\sum_{i,j}\eta^{ij}(\ell_i-\mu_i)\ell_j}$. As a result
$F(x)$ is still defined by (\ref{Fprod}), and the algebraic proofs remain
identical, apart from the use of the graded tensor products. The
analytic results are still satisfied without modification due to the
non degeneracy of $\eta$.

\indent

{\it Example: case of $osp(1|2)$}

\indent

In the case of $osp(1|2)$, with Cartan generator $h$, fermionic
step operators $e$ and $f$ and Cartan matrix $a^{\mbox{\tiny sym}}=(2)$, we have 
\begin{equation}
F_{12}(x)\ =\ \sum_{n=0}^{+\infty} e^n\otimes f^n \phi_n(x,h_{(2)}).
\end{equation}
The recursion relation satisfied by $(-1)^{n(n-1)/2}\phi_n(x,h)$ is
the same as the one found for $sl(2)$ provided that we change 
$q^2$ into $-q^2$ (except in the factors $(q-q^{-1})^n$ 
and $q^{-2h_{(2)}}$ which are left
unchanged) and $x^2$ into $-x^2$. 
Finally, 
\begin{equation}
F_{12}(x)\ =\ \sum_{n=0}^{+\infty}
\frac{(q-q^{-1})^n}{[n]_{-q^{2}}!}\, e^n\otimes f^n
\frac{(-1)^{n(n+1)/2}}{\prod_{\nu=1}^{n}\left(1+x^{-2}(-q^{2})^\nu 
 q^{-2h_{(2)}}\right)}.
\end{equation}
and, for $x$ sufficiently large,
\begin{equation}
F(x) = \prod_{k=0}^{+\infty} B_2(x)^k {\hat R}_{12}^{-1} B_2(x)^{-k}
 = \prod_{k=0}^{+\infty}\ \exp_{-q^{2}}\left( -(q-q^{-1})\ x^{-2k}
  q^{-2k} q^{-2kh_{(2)}}\ e\otimes f\right).
\end{equation}

\section{Conclusion}

In our work we have obtained a universal solution of the GNF equations for
finite dimensional Lie (super)algebras.
There are different paths along which our work can be pursued.

It would be very interesting to generalize our work to the case of affine
Lie algebra. A step in this direction has been achieved by C. Fr{\o}nsdal
\cite{Fr}.
This should shed some light on the construction of elliptic
solutions of GNF equations.

It is now certainly possible to generalize our work \cite{BR} to arbitrary
non compact quantum
complex group. In particular explicit formulas for $F(x)$ should allow us
to
understand the proof of the Plancherel formula in these cases. 

It has been shown in \cite{Ba} that in the $sl(2)$ case there exists an
element $g(x)\in \mathfrak{U}_q((sl(2))$ such that
$F(x)=\Delta(g(x))g_2(x)^{-1}g_1(xq^{\ell_{(2)}})^{-1}.$
It is conjectured that in the case of $sl(n)$ there exists  $g(x)\in
\mathfrak{U}_q((sl(n))$ such that
 \begin{equation}
F(x)=\Delta(g(x))F g_2(x)^{-1}g_1(xq^{\ell_{(2)}})^{-1}
\end{equation}
where $F$ satisfy the cocycle equation. 
The exact expression for $F(x)$ should be useful in the understanding of
this property.

Finally it would be very interesting to use our framework to derive
the eigenfunctions of Ruijsenaars--Schneider system (which are related 
to Macdonald's polynomials) using purely quantum
group  techniques.

\paragraph{Acknowledgments:}

We warmfully thank O.Babelon for numerous and fruitful discussions.

\bibliographystyle{unsrt}

\end{document}